\documentclass[12pt,twoside,leqno]{article}
\usepackage{latexsym}
\def\im{\mathop{\rm Im}\nolimits}
\def\dom{\mathop{\rm Dom}\nolimits}

\def\implies{\; \Longrightarrow \;}

\begin{document}
\begin{center}
{\bf COMBINATORIAL RESULTS FOR CERTAIN SEMIGROUPS OF ORDER-PRESERVING FULL CONTRACTION MAPPINGS OF A FINITE CHAIN}\\[4mm]
\textbf{A. D. Adeshola and A. Umar}\\

\end{center}

\newtheorem{theorem}{{\bf Theorem}}[section]
\newtheorem{prop}[theorem]{{\bf Proposition}}
\newtheorem{lemma}[theorem]{{\bf Lemma}}
\newtheorem{corollary}[theorem]{{\bf Corollary}}
\newtheorem{remark}[theorem]{{\bf Remark}}
\newtheorem{conj}[theorem]{{\bf Conjecture}}
\newcommand{\pf}{\smallskip\noindent {\em Proof.}\ \  }
\newcommand{\qed}{\hfill $\Box$\medskip}

\newcommand{\inv}{^{-1}}

\begin{abstract}
Let ${\cal T}_n$ be the  full symmetric semigroup on $X_n = \{1,
2, \ldots , n\}$ and let ${\cal OCT}_n$ and ${\cal ORCT}_n$ be its
subsemigroups of order-preserving and order-preserving or order-reversing
full contraction mappings of $X_n$, respectively. In this
paper we investigate the cardinalities of some equivalences on
${\cal OCT}_n$ and ${\cal ORCT}_n$ which lead naturally to obtaining
the orders of these subsemigroups. 
\footnote{This work was begun when the first named author was
visiting Sultan Qaboos University for a 3-month research visit in
Fall 2012.}\end{abstract}


\section{Introduction and Preliminaries}

Let $X_n = \{1, 2, \ldots , n\}$. A (partial) transformation $\alpha
:\dom \alpha \subseteq  X_n \rightarrow \im \alpha \subseteq X_n$ is
said to be {\em full} or {\em total} if $\dom \alpha = X_n$;
otherwise it is called {\em strictly} partial. The set of of full transformations of
$X_n$, denoted by ${\cal T}_n$, more commonly known as the {\em full transformation
semigroup} is also known as the {\em full symmetric semigroup or monoid} with
{\em composition of mappings} as the semigroup operation. We shall write $x\alpha$ for
the image of $x$ under $\alpha$ instead of $\alpha(x)$. The {\em height} of $\alpha$ is
denoted and defined by $h(\alpha)=\mid\mbox{Im}\,\alpha\mid$, the {\em right [left] waist}
of $\alpha$ is denoted and defined by $w^+(\alpha) = \max(\im \alpha)\, [w^-(\alpha) =
\min(\im \alpha)],$ \noindent the {\em fix} of $\alpha$ is denoted and defined by
$f(\alpha)= \mid F(\alpha)\mid$ where

$$F(\alpha) = \{x\in \dom \alpha: x\alpha = x\}.$$

It is also well-known that a partial transformation $\epsilon$ is
{\em idempotent} (${\epsilon}^2= \epsilon$) if and only if
$\im \epsilon = F(\epsilon)$. It is worth noting that to define the
left (right) waist of a transformation the base set $X_n$ must be totally
ordered.

\noindent A transformation $\alpha \in {\cal T}_n$ is said to be
{\em order-preserving (order-reversing)} if $(\forall x,y \in \
\dom \alpha)\ x \leq y \implies x\alpha \leq y\alpha \ (x\alpha
\geq y\alpha)$ and, a {\em contraction mapping} (or simply a {\em contraction}) if
($\forall x,y \in \dom \alpha) \mid x-y\mid \geq \mid x\alpha -y\alpha\mid$.
We shall denote by ${\cal OCT}_n$ and ${\cal ORCT}_n$, the semigroups
of order-preserving full contractions and of order-preserving or order-reversing
full contractions of an $n-$chain, respectively.

Recently, Zhao and Yang \cite{Zha} initiated the algebraic study of semigroups
of order-preserving partial contractions of an $n$-chain, where they referred to our
contractions as {\em compressions}. This paper investigates the combinatorial properties
of ${\cal OCT}_n$, ${\cal ORCT}_n$ and ${\cal ODCT}_n,$ and of their subsets of idempotents.

In this section we introduce basic terminologies and prove some preliminary results.
In Section 2 we obtain the cardinalities of various equivalence classes
defined on ${\cal OCT}_n$ and in Sections 3 and 4 we obtain the analogues of the results
in Section 2 for ${\cal ORCT}_n$ and ${\cal ODCT}_n$, respectively. These cardinalities lead to
formulae for the orders of ${\cal OCT}_n$, ${\cal ORCT}_n$ and ${\cal ODCT}_n$ as well
as new triangles of numbers that as a result of this work were recently recorded in
\cite{Slo}.

For standard concepts in semigroup and transformation semigroup
theory, see for example \cite{How3, Gan}. Let
\begin{eqnarray} \label{eqn1.1} {\cal OR}_n= \{\alpha \in {\cal T}_n:
(\forall x, y \in X_n)\, x\leq y \implies x\alpha\leq y\alpha\,\, \mbox{or}
\,\, x\alpha\geq y\alpha\}
\end{eqnarray}
\noindent be the subsemigroup of ${\cal T}_n$ consisting of all
order-preserving or order-reversing full transformations of $X_n$, and let
\begin{eqnarray} \label{eqn1.2} {\cal O}_n= \{\alpha \in {\cal T}_n:
(\forall x, y \in X_n)  x\leq y \implies x\alpha\leq y\alpha\}
\end{eqnarray}
\noindent be the subsemigroup of ${\cal T}_n$ consisting of all
order-preserving full transformations of $X_n$. Also let
\begin{eqnarray} \label{eqn1.3} {\cal CT}_n= \{\alpha \in {\cal T}_n:
(\forall x, y \in X_n)\, \mid x - y \mid \geq \mid x\alpha - y\alpha\mid\}
\end{eqnarray}
\noindent be the subsemigroup of ${\cal T}_n$ consisting of all
full contractions of $X_n$, and let
\begin{eqnarray} \label{eqn1.4} {\cal D}_n= \{\alpha \in {\cal T}_n:
(\forall x\in X_n)\, x\alpha \leq x\}
\end{eqnarray}
\noindent be the subsemigroup of ${\cal T}_n$ consisting of all
{\em order-decreasing} full transformations of $X_n$.

\noindent We have the following results

\begin{lemma} \label{lem1.1} Let $\alpha \in {\cal CT}_n$ be such that $f(\alpha)=m$.
Then $F(\alpha)=\{i, i+1, \ldots, i+m-1\}.$ Equivalently, $F(\alpha)$ is convex.
\end{lemma}

\pf Observe that it is sufficient to show that any point between two fixed points
(of $\alpha \in {\cal CT}_n$) must also be a fixed point. Let $x, y \in F(\alpha).$ Then
$x\alpha = x$ and  $y\alpha = y$. Suppose also without loss of generality $x < x' < y$
for some $x' \in X_n$. Note that if $x' =x'\alpha$, there is nothing to prove. Thus we consider
two cases: (i) $x'>x'\alpha$; \, (ii) $x'<x'\alpha$. In the former, we have
$$\mid y-x'\mid = \mid y\alpha-x'\mid = y\alpha - x'< y\alpha -x'\alpha=
\mid y\alpha - x'\alpha\mid$$
\noindent which implies that $\alpha$ is not a contraction. Hence we get a contradiction.
In the latter, we have
$$\mid x'-x\mid = \mid x' - x\alpha\mid = x' - x\alpha < x'\alpha -x\alpha=
\mid x'\alpha -x\alpha\mid$$
\noindent which implies that $\alpha$ is not a contraction. Hence we get a contradiction. Thus,
the proof is complete.\qed

\begin{lemma} \label{lem1.2} Let $\alpha \in {\cal CT}_n$ be such that $h(\alpha)=p$.
Then $\im \alpha=\{i, i+1, \ldots, i+p-1\}.$ Equivalently, $\im \alpha$ is convex.
\end{lemma}

\pf (By contradiction) Suppose that $\im \alpha$ is not convex. Then there exist
$x,z\in \im \alpha$ with $x<y<z$ for some $y\in X_n\setminus \im \alpha.$
Let $(y-1]$ and $[y+1)$ be the lower and upper saturations of $y-1$ and $y+1$,
respectively. Notice that $x\in (y-1]$ and $z\in [y+1).$ Moreover,
$(y-1]\alpha^{-1} \neq X_n \neq [y+1)\alpha^{-1}$ but $(y-1]\alpha^{-1} \cup [y+1)\alpha^{-1}= X_n.$
If  $(y-1]\alpha^{-1}$ is convex then since $(y-1]\alpha^{-1} \neq X_n$ there exist either
(i) $t\in (y-1]\alpha^{-1}$ and $t+1\in [y+1)\alpha^{-1}$; or
(ii) $t\in (y-1]\alpha^{-1}$ and $t-1\in [y+1)\alpha^{-1}.$

Case (i): it is clear that $t\alpha\leq y-1$ and $(t+1)\alpha \geq y+1$ so that
$$2\leq (t+1)\alpha - t\alpha = \mid (t+1)\alpha - t\alpha \mid \leq \mid (t+1) - t \mid = 1,$$
\noindent which is a contradiction.

Case (ii): it is clear that $t\alpha\leq y-1$ and $(t-1)\alpha \geq y+1$ so that
$$2\leq (t-1)\alpha - t\alpha = \mid (t-1)\alpha - t\alpha \mid \leq \mid (t-1) - t \mid = 1,$$
\noindent which is another contradiction.

On the other hand if $(y-1]\alpha^{-1}$ is not convex then there exists $t\in (y-1]\alpha^{-1}$
and either $t+1\in [y+1)\alpha^{-1}$ or $t-1\in [y+1)\alpha^{-1}.$ In the former, we see that
$t\alpha\leq y-1$ and $(t+1)\alpha\geq y+1.$ Thus,
$$2\leq (t+1)\alpha - t\alpha = \mid (t+1)\alpha - t\alpha \mid \leq \mid (t+1) - t \mid = 1,$$
\noindent which is a contradiction. In the latter, we see that
$t\alpha\leq y-1$ and $(t-1)\alpha\geq y+1.$ Thus,
$$2\leq (t-1)\alpha - t\alpha = \mid (t-1)\alpha - t\alpha \mid \leq \mid (t-1) - t \mid = 1,$$
\noindent which is another contradiction. Hence, the proof is complete.\qed

Next we state two important identities that will be needed later.

\begin{lemma} \label{lem1.3}(Vandemonde's Convolution Identity, \cite[(3a), p.8]{Rio}).
For all natural numbers m, n and p we have

$$\sum_{k=0}^{n}{{n\choose {m-k}}{p\choose k}}={n+p\choose m}.$$

\end{lemma}

\begin{lemma} \label{lem1.4}\cite[Lemma 1.3]{Lar3}
For all natural numbers $j$ and $a$ we have

$$\sum_{i=0}^{j-a}{j-i\choose a}={j+1 \choose a+1}.$$

\end{lemma}

\begin{lemma} \label{lem1.5}(\cite[(3b), p.8]{Rio}).
For all natural numbers m, n and p we have

$$\sum_{k=0}^{n}{{n-k\choose n-m}{p+k-1\choose k}}={n+p\choose m}.$$

\end{lemma}

\section{Order-preserving Full Contractions}

Enumerative problems of an essentially combinatorial nature arise
naturally in the study of semigroups of transformations. Many
numbers and triangle of numbers regarded as combinatorial gems like
the Stirling numbers \cite[pp. 42 \& 96]{How2}, the factorial
\cite{Uma1}, the Fibonacci number \cite {How1}, Catalan numbers
\cite{Gan}, Lah numbers \cite{Gan}, etc., have all featured in these
enumeration problems. For a nice survey article concerning combinatorial
problems in the partial transformation semigroup and some of its subsemigroups
we refer the reader to Umar \cite{Uma3}. These enumeration problems lead
to many numbers in Sloane's encyclopaedia of integer sequences \cite{Slo}
but there are also others that are not yet or have just been recorded in
\cite{Slo}.

\noindent Now recall the definitions of {\em height}, {\em (right) waist} and {\em fix}
of $\alpha\in {\cal T}_n$ stated earlier. From Umar \cite{Uma3}, we quote this result.

\begin{lemma} Let $X_n=\{1,2, \ldots,n\}$ and $P=\{p,m,k\}$, where for a given
$\alpha\in {\cal T}_n$, we set $p=h(\alpha), m=f(\alpha)$ and $k=w^+(\alpha)$.
Then we have the following:
\begin{enumerate}
    \item $n\geq k\geq p\geq m\geq 1$;
    \item $k=1 \implies p=1 \implies m\leq 1$.
 \end{enumerate}
\end{lemma}

As in Umar \cite{Uma3} let $S$ be a set
full transformations on $X_n$ and consider the combinatorial
functions:

$$F(n;p,m,k)= \mid\{\alpha \in S: \wedge (h(\alpha)=p, f(\alpha)= m, w^{+}(\alpha)=k)\}\mid,$$
$$F(n;p,m)= \mid\{\alpha \in S: \wedge (h(\alpha)=p, f(\alpha)= m\}\mid,$$
$$F(n;p,k)= \mid\{\alpha \in S: \wedge (h(\alpha)=p, w^{+}(\alpha)=k)\}\mid,$$
$$F(n;m,k)= \mid\{\alpha \in S: \wedge (f(\alpha)= m, w^{+}(\alpha)=k)\}\mid,$$
$$F(n;p)= \mid\{\alpha \in S: h(\alpha)=p\}\mid,$$
$$F(n;m)= \mid\{\alpha \in S: f(\alpha)= m\}\mid,$$
$$F(n;k)= \mid\{\alpha \in S: w^{+}(\alpha)=k\}\mid.$$

\noindent Observe that

$$F(n;a_1,a_2)= \sum_{a_3}{F(n;a_1,a_2,a_3)}, \,\, F(n;a_1)=
\sum_{a_2}{F(n;a_1,a_2)},$$

\noindent and

$$|S|= \sum_{a_1}{F(n;a_1)}$$
\noindent where $\{a_1,a_2,a_3\}= \{p,m,k\}.$

The following lemma is crucial to our investigation.

\begin{lemma} \label{lem2.2} Let
\begin{eqnarray*}a_{n,p} & = & \mid\{\alpha\in {\cal OCT}_n: (F(\alpha)=\{1\})\wedge (h(\alpha)=p)\}\mid\\
& = & \mid\{\alpha\in {\cal OCT}_n: (F(\alpha)=\{n\})\wedge (h(\alpha)=p)\}\mid.
\end{eqnarray*}
\noindent Then $a_{n,p} = {n-2\choose p-1}.$
\end{lemma}

\pf Let
$\alpha\in {\cal OCT}_n$ be such that  $F(\alpha)=\{1\}$ and  $h(\alpha)=p$. First observe that
$\im \alpha = \{1,2, \ldots, p\},$ by Lemma \ref{lem1.2}. Next, recall that the {\em blocks}
($y\alpha^{-1}$ for $y\in \im \alpha$) are convex by order-preservedness and so the number of
partitions of $X_n$ into $p$ convex blocks is the number of ways of inserting $p-1$ bars between
the $n-1$ spaces (between the points of the $n$-chain). However, notice that in this case the space
between 1 and 2 is not available so that 2 does not become a fixed point since $F(\alpha)=\{1\}$.
Thus, we have $n-2$ available places to insert $p-1$ bars, which can be done in ${n-2\choose p-1}$
ways, without new points being introduced. The result when $F(\alpha)=\{n\}$, follows by symmetry.\qed

The main result in this paper can now be stated.

\begin{prop}\label{prop1} Let $S={\cal OCT}_n = {\cal O}_n\cap {\cal CT}_n$. Then
$ F(n;p,m,k)=  {n-m-1\choose n-p-1}.$
\end{prop}

\pf Let $\alpha\in {\cal OCT}_n$ be such that  $h(\alpha)=p$, $f(\alpha)=m$ and $w^+(\alpha)=k$.
Then by Lemmas \ref{lem1.1} \& \ref{lem1.2} we see that

$F(\alpha)=\{i,i+1 \dots, i+m-1\}\subseteq \{k-p+1,k-p+2 \dots, k\}=\im \alpha,$

\noindent where $k-p+1\leq i\leq k-m+1.$ Next, observe that by order-preservedness we may decompose
$\alpha$ into $\beta =\alpha|_{\{1,2, \dots, i\}}, \, id_{\{i+1,i+2, \dots, i+m-2\}}$ and
$\beta'= \\ \alpha|_{\{i+m-1,i+m, \dots, n\}},$ where $F(\beta)=\{i\}, \,
\im \beta = \{k-p+1, \ldots, i\}, \, F(\beta')=\{i+m-1\},$ and $\im \beta' = \{i+m-1, \ldots, k\}.$ Now it
is not difficult to see that for $\beta$ there are $a_{i,i-k+p}$ possible maps, and for $\beta'$ there are
$a_{n-i-m+2,k-i-m+2}$ possible maps. Thus, multiplying the two numbers and taking
the sum of the product over $i$ from $k-p+1$ to $k-m+1$ and using Lemmas \ref{lem1.5} \& \ref{lem2.2} gives

 \begin{eqnarray*} F(n;p,k,m) & = & \sum_{i=k-p+1}^{k-m+1}a_{i,i-k+p}a_{n-i-m+2,k-i-m+2}\\
& = & \sum_{i=k-p+1}^{k-m+1}{i-2\choose i-k+p-1}{n-i-m\choose k-i-m+1} \\
& = & \sum_{i=k-p+1}^{k-m+1}{i-2\choose i-k+p-1}{n-i-m\choose n-k-1} \\
& = & {n-m-1\choose p-m} = {n-m-1\choose n-p-1},
\end{eqnarray*}
\noindent using the substitution $i-k+p-1=j$, for the last step.
\qed

\begin{corollary}\label{cor2.1} Let $S={\cal OCT}_n$. Then $F(n;p,m)=
(n-p+1){n-m-1\choose n-p-1},$ for $p\geq m\geq 1$.
\end{corollary}

\begin{corollary}\label{cor2.2} Let $S={\cal OCT}_n$. Then $F(n;p,k)=
{n-1\choose p-1}$, for $k\geq p\geq 1$.
\end{corollary}

\pf It follows directly from Proposition \ref{prop1} and Lemma \ref{lem1.4}.
\qed

\begin{corollary}\label{cor2.3} $S={\cal OCT}_n$. Then $F(n;m,k)=
\sum_{p=m}^{k}{n-m-1\choose n-p-1}$, for $k\geq m\geq 1$.
\end{corollary}

\pf It follows directly from Proposition \ref{prop1} and Lemma \ref{lem1.4}.
\qed

\begin{corollary}\label{cor2.4} Let $S={\cal OCT}_n$. Then $F(n;p)=
(n-p+1){n-1\choose p-1}$, for $p\geq 1$.
\end{corollary}

\pf It follows from either of the Corollaries \ref{cor2.1} \& \ref{cor2.2}.
\qed

\begin{corollary}\label{cor2.5} Let $S={\cal OCT}_n$. Then $F(n;k)=
\sum_{p=1}^{k}{n-1\choose p-1}$, for $k\geq 1$.
\end{corollary}

\pf It follows from either of the Corollaries \ref{cor2.2} \& \ref{cor2.3}.
\qed

For $0\leq i\leq n$, let

$$F(n;m_i)= \mid\{\alpha \in S: f(\alpha)= i\}\mid.$$

\begin{corollary}\label{cor2.6} Let $S={\cal OCT}_n$. Then $F(n;m_n)=1$
and $F(n;m)= \\(n-m+3)2^{n-m-2},$ for $n\geq 2$ and $n> m\geq 1$.
\end{corollary}

\pf It follows from either of the Corollaries \ref{cor2.1} \& \ref{cor2.3}.
\qed

\begin{corollary}\label{cor2.7} Let $S={\cal OCT}_n$. Then
$\mid S\mid=\mid{\cal OCT}_n\mid = (n+1)2^{n-2},$ for $n\geq 1$.
\end{corollary}

\pf It follows from any one of the Corollaries \ref{cor2.4}, \ref{cor2.5} \& \ref{cor2.6}.
\qed

\begin{corollary}\label{cor2.8} Let $S=E({\cal OCT}_n)$. Then $F(n;p,k)
= F(n;m,k) = 1,$ for $k\geq p = m\geq 1$.
\end{corollary}

\pf Since $F(\alpha) = \im \alpha$ for idempotents, it follows that $p=m$.
Hence the result follows from Proposition \ref{prop1}.
\qed

\begin{corollary}\label{cor2.9} Let  $S=E({\cal OCT}_n)$. Then $F(n;p)= F(n;m)=
n-p+1 = n-m+1$, for $p=m\geq 1$.
\end{corollary}

\begin{corollary}\label{cor2.10} Let  $S=E({\cal OCT}_n)$. Then $F(n;k)= k$,
for $k\geq 1$.
\end{corollary}

\begin{corollary}\label{cor2.11}
$\mid E({\cal OCT}_n)\mid = n(n+1)/2 = {n+1\choose 2},$ for $n\geq 1$.
\end{corollary}

\begin{remark} The triangle of numbers $F(n;m)$ has as a result of this work been
recorded in \cite{Slo}. However, $F(n;p)$, $F(n;k)$, $F(n;m_1)$ and $\mid {\cal OCT}_n\mid$
were recorded (in \cite{Slo}) as A104698, A008949, A045623 and A001792, respectively.
\end{remark}

\section{Order-preserving or Order-reversing Full Contractions}

\begin{remark} \label{rem1} For $h(\alpha)=p=1$ the concepts of
order-preserving and order-reversing coincide but distinct
otherwise. However, the map $\alpha \longmapsto \alpha h,$ where $xh=n-x+1$, for all $x$ in
$X_n$ is a bijection between the two sets for $p\geq 2$, see \cite[page 2, last paragraph]{Fer}.
\end{remark}

\begin{remark} Every idempotent is necessarily order-preserving.
Thus, there are no additional idempotents from reversing the order.
\end{remark}

The main result of this section is

\begin{prop}\label{prop2} Let $S={\cal ORCT}_n = {\cal OR}_n\cap {\cal CT}_n$. Then

$F(n;p,k)= \left\{
\begin{array}{ll}
2{n-1\choose p-1},\,&\,\,
\,p>1; \\
\,\,\,\,\,1,\,&\,\,\mbox{otherwise}.
\end{array}
\right.$

\end{prop}

\pf It follows from Corollary \ref{cor2.2} and Remark \ref{rem1}.
\qed

\begin{corollary}\label{cor3.1} Let $S={\cal ORCT}_n$. Then

$F(n;p)= \left\{
\begin{array}{ll}
2(n-p+1){n-1\choose p-1},\,&\,\,
\,p>1; \\
\,\,\,\,\,n,\,&\,\,\mbox{otherwise}.
\end{array}
\right.$
\end{corollary}

\pf It follows from Proposition \ref{prop2} and the fact that $p\leq k\leq n.$
\qed

\begin{corollary}\label{cor3.2} Let $S={\cal ORCT}_n$. Then $F(n;k)=
2\sum_{p=1}^{k}{n-1\choose p-1}-1$, for $k\geq 1$.
\end{corollary}

\pf It follows from Proposition \ref{prop2} and the fact that $1\leq p\leq k.$
\qed

\begin{corollary}\label{cor3.3} Let $S={\cal ORCT}_n$. Then
$\mid S\mid=\mid{\cal ORCT}_n\mid = (n+1)2^{n-1}-n,$ for $n\geq 1$.
\end{corollary}

\pf It follows from either of the Corollaries \ref{cor3.1} \& \ref{cor3.2}.
\qed

Let us denote the set of all order-reversing full contraction of $X_n$ by ${\cal ORCT}_n^*$
and let $$b(n,p)=\mid\{\alpha\in {\cal ORCT}_n^*:h(\alpha)=p\,\, \mbox{and}\, f(\alpha)=1\}\mid,$$
\noindent be the number of order-reversing full contractions of height $p$ and with exactly one
fixed point. Then we have

\begin{lemma}\label{lem3.4} For $n\geq p\geq 1$, we have $b(n,p)= (n-p+1)\sum_{i\geq 1}{n-2i\choose p-2i}$.
\end{lemma}

\begin{lemma}\label{lem3.5} For $n\geq p\geq 2$, we have
 \begin{itemize}
\item [(a)] $b(n,p)= (n-p+1){n-2\choose p-1} +b(n-2,p-2),\,\,\,b(n,1)=n, b(2,2)=0$,
\item [(b)] $(n-p)b(n,p)= (n-p)b(n-1,p-1)+(n-p+1)b(n-2,p-2),\\
b(2r+1,2r+1)=1, b(2r,2r)=0$.
\end{itemize}
\end{lemma}

Define a sequence $\{a_n\}$ by $$na_n=(n+2)a_{n-1}+2(n+1)a_{n-2},\,\,\, a_1=1, a_2=2.$$
Then we have

\begin{lemma}\label{lem3.6} Let $S={\cal ORCT}_n^*$. Then
 \begin{itemize}
\item[(a)] $F(n;m_0)=a_{n-1},\,\, n\geq 2;$
\item[(b)] $F(n,m_1)=a_n,\,\, n\geq 1;$
\item[(c)] $F(n;m)= 0,\,\, m\geq 2.$
\end{itemize}
\end{lemma}

\begin{lemma}\label{lem3.7} Let $S={\cal ORCT}_n$. Then $F(1;m_0)=0,\, F(n;m_n)=1$ \, and
 \begin{itemize}
\item[(a)] $F(n;m_0)=a_{n-1},\,\, n\geq 2;$
\item[(b)] $F(n,m_1)=a_n+(n+2)2^{n-3} - n,\,\, n\geq 1;$
\item[(c)] $F(n;m)= (n-m+3)2^{n-m-2},\,\, n>m\geq 2.$
\end{itemize}
\end{lemma}

\begin{remark} The triangles of numbers $F(n;p)$, $F(n;k)$ and $F(n;m)$; and the
sequences $F(n;m_0)$, $F(n;m_1)$ and $\mid {\cal ORCT}_n\mid$ have as a result of this work
just been recorded in \cite{Slo}.
\end{remark}

\section{Order-preserving and Order-decreasing Full Contractions}

The main result of this section is

\begin{prop}\label{prop3} Let $S={\cal ODCT}_n = {\cal D}_n\cap {\cal OCT}_n$. Then

$F(n;p,k,m)= \left\{
\begin{array}{ll}
{n-m-1\choose p-m},\,&\,
\,\,p=k; \\
\,\,0,\,&\,\,\mbox{otherwise}.
\end{array}
\right.$
\end{prop}

\pf Let $\alpha\in {\cal ODCT}_n$ be such that  $h(\alpha)=p$, $f(\alpha)=m$ and $w^+(\alpha)=k$.
Then by Lemmas \ref{lem1.1} \& \ref{lem1.2} we see that

$F(\alpha)=\{1, 2 \dots, m\}\subseteq \{1,2 \dots, p\}=\im \alpha.$

\noindent Moreover, $p=k$. Next, observe that by order-preservedness we may decompose
$\alpha$ into $id_{\{1,2, \dots,m-1\}}$ and $\beta'=\alpha|_{\{m,m+1, \dots, n\}},$ where
$F(\beta')=\{m\},$ and $\im \beta' = \{m,m+1 \ldots, p\}.$ Now it
is not difficult to see that there are
$a_{n-m+1,p-m+1}$ possible  $\beta'$s, by Lemma \ref{lem2.2}. Thus,

$F(n;p,k,m)= \left\{
\begin{array}{ll}
{n-m-1\choose p-m},\,&\,
\,p=k; \\
\,\,\,\,\,\,\,\,\,0,\,&\,\,\mbox{otherwise}.
\end{array}
\right.$ \qed

\begin{corollary}\label{cor4.1} Let $S={\cal ODCT}_n$. Then $F(n;p,m)=
{n-m-1\choose p-m},$ for $p\geq m\geq 1$.
\end{corollary}

\begin{corollary}\label{cor4.2} $S={\cal ODCT}_n$. Then $F(n;m,k)=
{n-m-1\choose k-m}$, for $k\geq m\geq 1$.
\end{corollary}

\begin{corollary}\label{cor4.3} Let $S={\cal ODCT}_n$. Then $F(n;p)=
{n-1\choose p-1}$, for $p\geq 1$.
\end{corollary}

\begin{corollary}\label{cor4.4} Let $S={\cal ODCT}_n$. Then $F(n;k)=
{n-1\choose k-1}$, for $k\geq 1$.
\end{corollary}

\begin{corollary}\label{cor4.5} Let $S={\cal ODCT}_n$. Then $F(n;m_n)=1$
and  $F(n;m)= 2^{n-m-1},$ for $n > m\geq 1$.
\end{corollary}

\begin{corollary}\label{cor4.6} Let $S={\cal ODCT}_n$. Then
$\mid S\mid=\mid{\cal ODCT}_n\mid = 2^{n-1},$ for $n\geq 1$.
\end{corollary}

\begin{corollary}\label{cor4.7} Let $S=E({\cal ODCT}_n)$. Then

$F(n;p,m)= \left\{
\begin{array}{ll}
1,\,&\,\,p=m; \\
0,\,&\,\,\mbox{otherwise}.
\end{array}
\right.$
\end{corollary}

\pf Since $F(\alpha) = \im \alpha$ for idempotents, it follows that $p=m$.
Hence the result follows from Proposition \ref{prop1}.
\qed

\begin{corollary}\label{cor4.8} Let  $S=E({\cal ODCT}_n)$. Then $F(n;p)= F(n;m)= F(n;k) =1$, for $p,k,m\geq 1$.
\end{corollary}

\begin{corollary}\label{cor4.9}
$\mid E({\cal ODCT}_n)\mid = n,$ for $n\geq 1$.
\end{corollary}

\noindent {\bf Acknowledgements}. The first named author would like
to thank Kwara State University, Malete and TET Fund for financial support, and Sultan Qaboos
University for hospitality.

\vspace{1cm}

\begin{center}
A.\ D.\ Adeshola\\
Department of Statistics and Mathematical Sciences\\
Kwara State University, Malete \\
P. M. B. 1530, Ilorin\\
Nigeria.\\
E-mail:{\tt adeshola.dauda@kwasu.edu.ng}

\end{center}

\begin{center}
 A.\ Umar\\
Department of Mathematics and Statistics\\
Sultan Qaboos University \\
Al-Khod, PC 123 -- OMAN\\
E-mail:{\tt aumarh@squ.edu.om}

\end{center}

\end{document}